\documentclass[11pt,reqno]{amsart}
\usepackage{amsmath}
\usepackage{amsxtra}
\usepackage{amscd}
\usepackage{amsthm}
\usepackage{amsfonts}
\usepackage{amssymb}
\usepackage{eucal}

\newcommand{\nn}{\nonumber}
\newcommand{\bea}{\begin{eqnarray}}
\newcommand{\ena}{\end{eqnarray}}
\newcommand{\be}{\begin{eqnarray*}}
\newcommand{\en}{\end{eqnarray*}}

\newcommand{\chib}{{\chi}^{W}}
\newcommand{\chiC}{{\chi}^{\mathcal{C}}}

\newcommand{\slrh}{\widehat{\mathfrak{sl}}_r}
\newcommand{\slth}{\widehat{\mathfrak{sl}}_2}
\newcommand{\slkh}{\widehat{\mathfrak{sl}}_k}
\newcommand{\hh}{\widehat{\mathfrak{h}}}
\newcommand{\Z}{{\mathbb Z}}  
\newcommand{\C}{{\mathbb C}}

\newcommand{\A}{\mathcal{A}}

\newtheorem{thm}{Theorem}[section]
\newtheorem{prop}[thm]{Proposition}
\newtheorem{lem}[thm]{Lemma}
\newtheorem{cor}[thm]{Corollary}

\numberwithin{equation}{section}

\begin{document}
\pagestyle{myheadings}
\markboth{Feigin, Jimbo, Miwa and Mukhin}
{Symmetric polynomials}

\title[Symmetric polynomials]
{Symmetric polynomials vanishing on the shifted diagonals 
and Macdonald polynomials
\\
}
\author{B. Feigin, M. Jimbo, T. Miwa and E. Mukhin}
\address{BF: Landau institute for Theoretical Physics, Chernogolovka,
142432, Russia}\email{feigin@feigin.mccme.ru}  
\address{MJ: Graduate School of Mathematical Sciences, University of
Tokyo, Tokyo 153-8914, Japan}\email{jimbomic@ms.u-tokyo.ac.jp}
\address{TM: Division of Mathematics, Graduate School of Science , 
Kyoto University, Kyoto 606-8502 Japan}\email{tetsuji@kusm.kyoto-u.ac.jp}
\address{EM: Department of Mathematics, 
Indiana University-Purdue University-Indianapolis, 
402 N.Blackford St., LD 270, 
Indianapolis, IN 46202}
\email{mukhin@math.iupui.edu}

\date{\today}

\setcounter{footnote}{0}\renewcommand{\thefootnote}{\arabic{footnote}}

\begin{abstract}
{}For each pair ($k, r$) of positive integers with $r\ge 2$,  
we consider an ideal $I^{(k,r)}_n$ of the ring of symmetric 
polynomials.  
The ideal $I_n^{(k,r)}$ has a basis consisting of Macdonald polynomials 
$P_\lambda(x_1,\cdots,x_n;q,t)$ at $t^{k+1}q^{r-1}=1$,  
and is a deformed version of the one studied earlier 
in the context of Jack polynomials. 
In this paper we give a characterization of 
$I^{(k,r)}_n$ in terms of explicit zero conditions on the 
$k$-codimensional shifted diagonals of the form 
$x_{2}=tq^{s_1}x_1,\cdots,x_{k+1}=tq^{s_k}x_k$. 

The ideal $I^{(k,r)}_n$ may be viewed as a deformation of 
the space of correlation functions of an abelian current of 
the affine Lie algebra $\slrh$. 
We give a brief discussion about this connection. 
\end{abstract}
\maketitle

\vskip20mm
\newpage

\setcounter{footnote}{0}
\renewcommand{\thefootnote}{\arabic{footnote})}
\renewcommand{\arraystretch}{1.2}

\setcounter{section}{0}
\setcounter{equation}{0}
\section{Introduction}
Recall that integrable representations of Kac-Moody Lie algebras 
can be characterized in terms of vanishing ideals. 
Let us formulate this fact for the current algebra $\slrh$. 
Denote by $e_{ij}(z)$ ($i\neq j$)  
the currents corresponding to 
root vectors $e_{ij}$ of $\mathfrak{sl}_r$. 
A representation of $\slrh$ in category $\mathcal{O}$ 
is integrable if and only if 
$e_{ij}(z)^{k+1}$ acts as $0$ for some $i\neq j$, where 
$k$ is a non-negative integer called the level. 
In this case $e_{ij}(z)^{k+1}=0$ holds for all $i\neq j$. 
Actually there are more relations which follow from these.
Consider the commutative currents 
$e_{12}(z),e_{13}(z),\cdots,e_{1r}(z)$. 
On level $k$ integrable representations they satisfy the relations 
\bea
&&
e_{12}(z)^{\nu_1}e_{13}(z)^{\nu_2}\cdots e_{1r}(z)^{\nu_{r-1}}=0
\label{arel}
\\
&&\quad\quad \mbox
{for all $\nu_1,\cdots,\nu_{r-1}\in\Z_{\ge 0}$
such that $\nu_1+\cdots +\nu_{r-1}=k+1$}.
\nn
\ena
Let us combine them into a single abelian current 
\bea
e(z)&=&e_{12}(z^{r-1})+z^{-1}e_{13}(z^{r-1})+\cdots+z^{-r+2}e_{1r}(z^{r-1}) 
\label{e-slr}\\
&=&\sum_{i\in\Z} e_i z^i. 
\nn
\ena
If we fix a primitive root of unity $\tau$ of order $r-1$,   
then the relations \eqref{arel} can be 
rewritten in terms of $e(z)$ as follows. 

\bea
&&e(z)^{\nu_1}e(\tau z)^{\nu_2}\cdots e(\tau^{r-2}z)^{\nu_{r-1}}=0
\nn\\
&&\qquad 
(\nu_1,\cdots,\nu_{r-1}\in\Z_{\ge 0}, 
\nu_1+\cdots+\nu_{r-1}=k+1). 
\label{erel}
\ena
Using \eqref{erel}, one can obtain 
a monomial basis of integrable representations \cite{P}. 
{}For simplicity, 
let us consider the vacuum representation $L$ of level $k$
with highest weight vector $v$ such that 
$e_i v=0$ for $i\le 0$.
We call $W=\C[e_1,e_2,\cdots]v\subset L$ 
the principal subspace of $L$. 
Then the following set of monomials constitutes a basis of $W$: 
\bea
&&e_1^{a_1}e_2^{a_2}\cdots e_l^{a_l}\quad (l\in \Z_{\ge 0}), 
\label{monb}\\
&&\mbox{where $a_i\in \Z_{\ge 0}$,  
$a_i+a_{i+1}+\cdots+a_{i+r-1}\le k$ for all $i\ge 1$}.
\nn
\ena

We present here a deformation of the relations \eqref{erel} 
which preserves the structure of the monomial basis \eqref{monb}. 
Namely, let $q,t$ be complex numbers such that 
$t^{k+1}q^{r-1}=1$. 
Consider an abelian current $e(z)=\sum_{i\in\Z}e_iz^i$ 
satisfying the relations  
\bea
&&e(z)e(tq^{s_1}z)\cdots e(t^kq^{s_k}z)=0\quad 
\label{qt}\\&&\quad \qquad
\mbox{for all $0\le s_1\le \cdots\le  s_{k}\le r-2$}.
\nn
\ena
{}For $t=1$ and $q=\tau$ we get back to the relations \eqref{erel}. 
If we replace $e(z)$ by $\sum_{i\ge 1}e_iz^i$ 
in \eqref{qt}, then the Fourier coefficients of the left hand side 
are well defined elements of $\C[\{e_i\}_{i\ge 1}]$. 
Let $\mathcal{J}(q,t)$ be the ideal generated by them. 
In this paper we prove that the quotient 
$\C[\{e_i\}_{i\ge 1}]/\mathcal{J}(q,t)$ 
has the same set \eqref{monb} 
as a monomial basis, 
if $q,t$ are `generic' (see \eqref{qtu}, \eqref{qatb} below 
for the precise condition). 
{}From this fact it follows that 
the current $e(z)$ of $\slrh$ acting on $L$ 
can be deformed to a current satisfying the relations \eqref{qt}. 
We also give an analogous result for 
integrable representations of $\slrh$ other than the vacuum module. 

{}From a slightly different point of view, our result can be described as follows. 
Let $D$ be an element of the Weyl group of $\slrh$ 
such that $D e_i D^{-1}=e_{i+r}$. 
Then, on level $k$ integrable representations of $\slrh$, 
we have an action of the algebra 
$E_{k,r}=\C[D,D^{-1}]\ltimes\C[\{e_i\}_{i\in \Z}]/
\tilde{\mathcal{J}}(\tau,1)$,   
where 
$\tilde{\mathcal{J}}(q,t)$ 
denotes the ideal of $\C[\{e_i\}_{i\in \Z}]$ 
generated by the Fourier coefficients of the left hand side of
\eqref{qt}. (To be precise a completion is necessary, 
but we do not discuss such details here.) 
It is possible to show that each level $k$ 
integrable irreducible representation of $\slrh$ 
remains irreducible upon restriction to $E_{k,r}$. 
(Actually $E_{k,r}$ 
has more irreducible representations than 
$\slrh$, but their classification is not known. )
Our result shows that for generic $q,t$ the algebra 
$E_{k,r}(q,t) 
=\C[D,D^{-1}]\ltimes\C[\{e_i\}_{i\in \Z}]/\tilde{\mathcal{J}}(q,t)$ 
has irreducible representations which are `deformations' of 
representations of $\slrh$. 

The relations of the type \eqref{qt} can be put in more general
context as follows.
{}Fix a function $\lambda(x,y)$ in two variables. 
{}Following \cite{OF} we define two associative algebras 
$S$ and $A$. 
The algebra $S$ is a graded algebra $S=\oplus_{n\ge 0}S_n$, 
each graded component $S_n$ being the space 
of symmetric functions in $n$ variables. 
{}For $F\in S_m$ and $G\in S_n$,  
the product $F*G\in S_{m+n}$ is defined by the formula
\bea
&&F*G(x_1,\cdots,x_{m+n})=
\nn\\
&&\mathop{\rm Sym}\Bigl(
F(x_1,\cdots,x_m)G(x_{m+1},\cdots,x_{m+n})
\prod_{1\le i\le m\atop m+1\le j\le m+n}
\lambda(x_i,x_j)\Bigr).
\label{star}
\ena
In \eqref{star}, the symbol $\mathop{{\rm Sym}}$ stands 
for the symmetrization. 
The algebra 
$A=\oplus_{n\ge 0}A_n$ 
is defined similarly, 
where $A_n$ 
consists of anti-symmetric functions 
in $n$ variables and symmetrization in \eqref{star} is replaced
by anti-symmetrization. 
Let $K\subset \C\times \C$ be the set of zeroes of the 
function $\lambda(x,y)$.
 
We say that $F\in S_n$ (or $F\in A_n$) 
satisfies the 
{\it wheel condition} if the following holds:
\bea
\mbox{$F=0$ whenever
$(x_1,x_2),\cdots,(x_{l},x_{l+1}),(x_{l+1},x_1)\in K$ for some
$1\leq l\leq n$}.
\label{wh}
\ena
If $F,G$ satisfy the wheel condition, then so does $F*G$. 
Therefore, functions satisfying the wheel condition \eqref{wh}
constitute a subalgebra $S^w\subset S$ (resp. $A^w\subset A$).
Note that the subalgebra generated by $S_1$ 
(resp. $A_1$) 
in $S$ (resp. $A$)
is contained in $S^w$ (resp. $A^w$).

Set 
\be
\lambda(x,y)=
\frac{(x-t_1y)\cdots(x-t_sy)}{(x-y)^s}
\en
where $\mathcal{S}=\{t_1,\cdots,t_s\}$ is a set of 
non-zero complex numbers.
We will refer to $\mathcal{S}$ as the {\it wheel set}.

Let $B(\mathcal{S})$ consist of rational functions
of the form 
\be
F(x_1,\cdots,x_n)=\frac{f(x_1,\cdots,x_n)}
{\prod_{1\le i<j\le n}(x_i-x_j)^{s-1}}, 
\en
where $f(x_1,\cdots,x_n)$ is a symmetric Laurent polynomial satisfying
\bea
&&\mbox{$f(x_1,\cdots,x_n)=0$  }
\nn\\
&&\qquad 
\mbox{if 
$\displaystyle{\frac{x_2}{x_1}=t_{i_1},
\cdots,\frac{x_{l+1}}{x_{l}}=t_{i_{l}},
\frac{x_1}{x_{l+1}}=t_{i_{l+1}}}$ 
for some $i_1,\cdots,i_{l+1}$
}.
\label{wh1}
\ena
Then $B(\mathcal{S})$ is a subalgebra of $S^w$
if $s$ is odd and of $A^w$ if $s$ is even. 

Note that the wheel condition is non-trivial 
only when the parameters $t_1,\cdots,t_s$  
satisfy 
\be
t_1^{\kappa_1}t_2^{\kappa_2}\cdots t_s^{\kappa_s}=1
\en
for some integers $\kappa_1,\cdots,\kappa_s$. 
We call such an equation a {\it resonance condition}.

Returning to the abelian current $e(z)$ and \eqref{qt},
we take $s=r$, $l=k$, and the wheel set
\bea
\mathcal{S}=\{t,tq,\cdots,tq^{r-1}\}
\label{S1}
\ena
where $t^{k+1}q^{r-1}=1$ is assumed.  
Then the relation \eqref{qt} coincides 
with \eqref{wh} in the following sense.

Suppose we have a representation $W$ 
of the current $e(z)$ satisfying \eqref{qt}.  
(For the relations \eqref{qt} to make sense, 
we consider only such representations that
for any vector $v\in W$ 
there is an integer $N$ satisfying $e_iv=0$ for $i<N$.) 
Then the matrix elements 
\bea
f(z_1,\cdots,z_n)=
\langle v^{\vee},e(z_1)\cdots e(z_n) v\rangle,
\label{G}
\ena
where $v\in W$ and $v^{\vee}\in W^*$, are Laurent polynomials 
satisfying the wheel condition \eqref{S1}.

Therefore the study of the representation $W$ is closely related to 
that of the algebra $B(\mathcal{S})$.
Now let $W=\C[e_1,e_2,\cdots]/\mathcal{J}(q,t)$, 
$v=1\bmod \mathcal{J}(q,t)$ and $e_iv=0$ ($i\le 0$). 
\footnote{In the main text, we shift the index of $e_i$ by one
and consider a quotient space of $\C[e_0,e_1,\cdots]$.} 
Then the dual space $W^*$ can be identified with the 
space $J^{(k,r)}=\oplus_{n\geq0}J^{(k,r)}_n$ of all symmetric polynomials 
satisfying the wheel condition relative to \eqref{S1}.
Here $J^{(k,r)}_n$ denotes the subspace of $n$ variables.

We use the theory of symmetric functions to find a basis in $J^{(k,r)}$.
Namely, let $I^{(k,r)}_n$ be spanned by
Macdonald polynomials $P_\lambda(x;q,t)$, where 
$q,t$ satisfy $t^{k+1}q^{r-1}=1$, 
and $\lambda$ ranges over a set of ($k,r,n$)-admissible (see
\eqref{adm part} for the definition) partitions.
Our main result is that 
for `generic' $q,t$ we have $I^{(k,r)}_n=J^{(k,r)}_n$.
In other words, the above Macdonald polynomials 
constitute a basis of $J^{(k,r)}_n$.
 
We do not understand well the reason why 
the particular choice \eqref{S1} of the wheel set 
is exactly what we need 
to deform the integrability condition \eqref{arel} for $\slrh$.

This paper is organized as follows. 
In Section 2, we review known facts about Macdonald polynomials.  
We discuss briefly their regularity properties
when the parameters $q,t$ satisfy the relation $q^at^b=1$
with some $a,b\in\Z_{\ge 1}$,  
following the work \cite{FJMM} on Jack polynomials.  
Our main result is stated as Theorem \ref{thm:1}.
Section 3 is devoted to its proof. 
In Section 4 we give a monomial base for the 
analogs of (non-vacuum) integrable representations of $\slrh$. 
In Section 5, we discuss an expected 
link between the present work and 
representations of $W_k$ algebras.

\setcounter{section}{1}
\setcounter{equation}{0}
\section{Macdonald polynomials and wheel condition}\label{sec:2}

\subsection{Preliminaries}
\label{subsec:2.1}
In this subsection, we review basic facts about 
the Macdonald polynomials which we use in the text. 
Our basic reference is Macdonald's book \cite{M}. 

Let $n$ be a non-negative integer. 
The Macdonald operators $\{D_n^r\}_{0\le r\le n}$ 
are mutually commuting $q$-difference operators acting on 
the ring of symmetric polynomials $\C(q,t)[x_1,\cdots,x_n]^{S_n}$, 
where $S_n$ stands for the symmetric group on $n$ letters. 
Explicitly they are given by the formula 
\bea
&&D_n^r=\sum_{|I|=r}A_I(x;t)T_{I},
\label{Dr}
\ena
where $I\subset\{1,\cdots,n\}$ runs over subsets of cardinality $r$,  
\be
&&A_I(x;t)=t^{r(r-1)/2}\prod_{i\in I \atop j\not\in I}
\frac{tx_i-x_j}{x_i-x_j}, \\
&&T_{I}=\prod_{i\in I}T_{q,x_i},
\en
and 
$(T_{q,x_i}f)(x_1,\cdots,x_n)=f(x_1,\cdots,qx_i,\cdots,x_n)$.
Let $D_n(X;q,t)=\sum_{r=0}^{n}D_n^rX^r$ be their generating function.

Denote by $\pi_n$ the set of partitions 
$\lambda=(\lambda_1,\ldots,\lambda_n)$,  
where $\lambda_i$ are non-negative integers 
satisfying $\lambda_i\geq\lambda_{i+1}$ ($1\le i\le n-1$).
The Macdonald polynomials $\{P_\lambda\}_{\lambda\in\pi_n}$ 
constitute a unique $\C(q,t)$-basis $\C(q,t)[x_1,\cdots,x_n]^{S_n}$ 
characterized by the following defining properties.
\bea
&&D_n(X;q,t)P_\lambda=\prod_{i=1}^n(1+Xq^{\lambda_i}t^{n-i})\cdot P_\lambda,
\label{diag}
\\
&&
P_\lambda=m_\lambda+\sum_{\mu<\lambda}u_{\lambda\mu}m_\mu 
\qquad (u_{\lambda\mu}\in\C(q,t)).  
\label{tri}
\ena
In the second line, 
$m_\lambda=\sum_{\alpha\in S_n\lambda}
x_1^{\alpha_1}\cdots x_n^{\alpha_n}$  
stands for the monomial symmetric function.
The dominance ordering $\mu<\lambda$ on $\pi_n$ is defined by 
$\mu\neq \lambda$ and 
$\mu_1+\cdots+\mu_i\le \lambda_1+\cdots+\lambda_i$ for all $i=1,\cdots,n$. 

Along with \eqref{Dr}, 
we also consider the following operators \cite{L}.  
{}For $m\ge 0$, set 
\bea
E_m=\sum_{i=1}^nx_i^m A_i(x;t)\frac{\partial}{\partial_q x_i},
\label{Em}
\ena
where $A_i(x;t)=A_{\{i\}}(x;t)$ and 
\be
&&\frac{\partial}{\partial_q x_i}=\frac{1}{(q-1)x_i}\left(T_{q,x_i}-1\right).
\en
{}For $m=1$, \eqref{Em} is related to $D_n^1$ by 
$(q-1)E_1=D_n^1-(1-t^n)/(1-t)$. 
Set further $e_1=\sum_{i=1}^nx_i$. 
Then we have 
\bea
&&
e_1P_\lambda=
\sum_{j=1}^{\ell(\lambda)+1}
\psi_{\lambda^{(j)}/\lambda}'P_{\lambda^{(j)}},
\label{e1}
\\
&&E_0P_\lambda=
\sum_{j=1}^{\ell(\lambda)}
\psi_{\lambda/\lambda_{(j)}}''P_{\lambda_{(j)}},
\label{E0}\\
&&E_2P_\lambda=\frac{t^{n-1}}{1-q}\sum_{j=1}^{\ell(\lambda)+1}
(1-q^{\lambda_j}t^{1-j})\psi_{\lambda^{(j)}/\lambda}'
P_{\lambda^{(j)}}.
\label{E2}
\ena
Here $\lambda^{(j)}$ (resp. $\lambda_{(j)}$) 
denotes the partition obtained 
by adding one node to 
(resp. removing one node from) the $j$-th row of $\lambda$,  
and $\ell(\lambda)=\max\{j\mid \lambda_j>0\}$ signifies the length of $\lambda$. 
The coefficients are given by 
\be
&&\psi_{\lambda^{(j)}/\lambda}'=
\prod_{i=1}^{j-1}
\frac{1-q^{\lambda_i-\lambda_j-1}t^{j-i+1}}{1-q^{\lambda_i-\lambda_j}t^{j-i}}
\frac{1-q^{\lambda_i-\lambda_j}t^{j-i-1}}{1-q^{\lambda_i-\lambda_j-1}t^{j-i}},
\\
&&\psi_{\lambda/\lambda_{(j)}}''=
\frac{1-t^{n-j}q^{\lambda_j}}{1-q} 
\prod_{i=j+1}^n
\frac{1-q^{\lambda_j-\lambda_i-1}t^{i-j+1}}{1-q^{\lambda_j-\lambda_i}t^{i-j}}
\frac{1-q^{\lambda_j-\lambda_i}t^{i-j-1}}{1-q^{\lambda_j-\lambda_i-1}t^{i-j}}.
\en
When $\lambda_j=\lambda_{j-1}$ and thus $\lambda^{(j)}$ is not defined, 
the corresponding term is absent in the right hand sides 
of \eqref{e1},\eqref{E2}, because $\psi_{\lambda^{(j)}/\lambda}'=0$. 
The same remark applies to \eqref{E0}. 
Eq.\eqref{e1} is a special case of the Pieri formula (\cite{M}, eq.(6.24)),
while \eqref{E0},\eqref{E2} are due to  \cite{L}. 

\subsection{Regularity of Macdonald polynomials}\label{subsec:2.2}
The coefficients $u_{\lambda\mu}$ in \eqref{tri}
are rational functions of $q$ and $t$. 
{}For a partition $\lambda\in\pi_n$, denote by
$\lambda'$ the conjugate partition. 
It is known \cite{Int} that if we set 
\be
c_\lambda(q,t)=\prod_{(i,j)\in\lambda}(1-q^{\lambda_i-j}t^{\lambda'_j-i+1}), 
\en
then $c_\lambda P_\lambda$ is a polynomial in $q,t$.
In particular, all possible poles of $u_{\lambda\mu}$ are of the form 
\bea
q^at^b=1\qquad (a,b\in\Z,~~ a\ge 0,b>0).  
\label{pole}
\ena
{}For given $n$ and $q,t$ satisfying \eqref{pole},  
it is natural to ask which $P_\lambda$ remain well defined. 
We have studied a sufficient condition in \cite{FJMM} 
in the limit $t=q^\beta$, $q\rightarrow 1$, 
where Macdonald polynomials reduce to Jack polynomials. 
As we show below, 
the results of \cite{FJMM} have straightforward 
extensions to the setting of Macdonald polynomials. 

Throughout this paper, 
we fix integers $k,r$ where $k\geq 1$ and $r\ge 2$. 
As opposed to \cite{FJMM}, we do not assume that $k+1,r-1$ are coprime. 
Let $m$ be the greatest common divisor of $k+1$ and $r-1$,   
and let $\omega$ be a primitive $m$-th root of unity. 
Let further $\omega_1\in\C$ be such that $\omega_1^{(r-1)/m}=\omega$. 
{}For an indeterminate $u$ we consider the specialization 
\bea
t=u^{\frac{r-1}{m}},
\quad 
q=\omega_1 u^{-\frac{k+1}{m}},
\label{qtu}
\ena
so that $t^{\frac{k+1}{m}}q^{\frac{r-1}{m}}=\omega$. 
{}For integers $a,b\in \Z$, we have then 
\bea
\mbox{$q^at^b=1$ if and only if 
$a=(r-1)s$, $b=(k+1)s$ for some $s\in\Z$}.  
\label{qatb}
\ena

As in \cite{FJMM}, we say that a partition $\lambda\in\pi_n$ is 
{\it $(k,r,n)$-admissible} if 
\bea\label{adm part}
\lambda_i-\lambda_{i+k}\ge r\qquad (i=1,\cdots,n-k). 
\ena
The following two Lemmas can be verified by noting \eqref{qatb} and 
repeating the working of Lemma 2.1--2.3 in \cite{FJMM}.

\begin{lem}\label{lem:2.1}
Suppose $1\le i<j\le n$ and $\lambda\in\pi_n$ is $(k,r,n)$-admissible. 
Then 
\be
&& q^{\lambda_i-\lambda_j}t^{j-i}\neq 1,\\
&& q^{\lambda_i-\lambda_j-1}t^{j-i+1}\neq 1,\\
&&q^{\lambda_i-\lambda_j-1}t^{j-i}\neq 1.\\
\en
If in addition $\lambda_j<\lambda_{j-1}$, then 
\be
&&q^{\lambda_i-\lambda_j}t^{j-i-1}\neq 1.
\en
\end{lem}

\begin{lem}\label{lem:2.11}
If $\lambda$ is ($k,r,n$)-admissible, then 
$\psi'_{\lambda^{(j)}/\lambda}$ is well defined. 
It is zero if and only if $\lambda_{j-1}=\lambda_j$.
\end{lem}

\begin{lem}\label{lem:2.2}
Assume either $\lambda$ is $(k,r,n)$-admissible, or else
$\lambda$ is obtained from a $(k,r,n)$-admissible 
partition by adding or removing one node. 
Then 
$P_\lambda$ has no pole 
at $(t,q)=(u^{\frac{r-1}{m}},\omega_1 u^{-\frac{k+1}{m}})$. 
\end{lem}
\begin{proof}
In the context of Jack polynomials, 
an analogous statement is given as Proposition 2.6 in \cite{FJMM}.
The same proof applies by using \eqref{qatb}, 
Lemma \ref{lem:2.1} and the formula \eqref{diag} 
for the eigenvalue of the Macdonald operators. 
We omit further details.
\end{proof}

{\it In the rest of this paper, we 
fix the specialization of $t,q$ as in \eqref{qtu}. }

\subsection{Statement of the result}
\label{subsec:2.3}

We set $\Lambda_n=K[x_1,\cdots,x_n]^{S_n}$,  
where the ground field is $K=\C(u)$. 
Define a subspace $I^{(k,r)}_n$ of $\Lambda_n$ by 
\be
I^{(k,r)}_n=\mbox{span}_{K}
\{P_\lambda(x_1,\cdots,x_n;q,t)\mid
\mbox{$\lambda$ is ($k,r,n$)-admissible}\}.
\en
Our goal is to characterize this space in terms of the wheel 
condition. 

Consider the subspace $J^{(k,r)}_n\subset \Lambda_n$ 
of all symmetric polynomials $f(x_1,\cdots,x_n)$ 
satisfying the wheel condition \eqref{wh1}
relative to the wheel set $\mathcal{S}=\{t,tq,\cdots,tq^{r-1}\}$. 
Equivalently, $f\in \Lambda_n$ belongs to
$J^{(k,r)}_n$ if and only if 
\bea
&&f=0~~~\mbox{ if }~~~~
x_i=t^{i-1}q^{s_1+\cdots+s_{i-1}}x_1
~~(2\le i\le k+1)
\nn\\&&
\mbox{for all $s_1,\cdots,s_{k+1}\in\Z_{\ge 0}$ 
satisfying $s_1+\cdots+s_{k+1}=r-1$.}
\label{zero}
\ena
Here we require the vanishing of $f$ on the wheel of length $k+1$.
Since we have the resonance of the form (\ref{qatb}), we have wheels of larger
length of a multiple of $k+1$. However, the vanishing of $f$ for
such a wheel follows from (\ref{zero}) because the larger wheel
necessarily contains a wheel of length $k+1$.

The following is our main result. 
\begin{thm}\label{thm:1} {}For all $n\ge 0$ we have an equality of 
ideals of symmetric polynomials 
\bea
I^{(k,r)}_n=J^{(k,r)}_n.
\label{main}
\ena
Moreover these ideals are stable under the action of 
$D(X;q,t)$ and $E_m$ ($m\ge 0$). 
\end{thm}
{}For $r=2$, the condition \eqref{zero} simplifies to 
\be
\mbox{$f=0$ if $x_j=t^{j-1}x_1$ for $j=1,\cdots,k+1$}.
\en
Theorem \ref{thm:1} in this case was stated in \cite{FJMM}.
\footnote{
In \cite{FJMM}, 5 lines above Theorem 4.4,
the condition ` $x_j=t^{j-1}$' should read ` $x_j=t^{j-1}x_1$'.
}

The next Section is devoted to the proof of Theorem \ref{thm:1}. 

\section{Proof of Theorem \ref{thm:1}}\label{sec:3}
\subsection{Stability by Macdonald type operators}
\label{subsec:3.2}

{}From the proof of Proposition 3.4--3.6 in \cite{FJMM} 
we see that, 
if $\lambda$ is $(k,r,n)$-admissible, then
the formulas \eqref{e1}--\eqref{E2} remain valid
if 
only admissible partitions are retained in the right hand side. 
In particular, 
the space $I^{(k,r)}_n$ is invariant under the action of 
the operators $D_n(X;q,t)$, $E_m$  ($m=0,1,2$) 
and multiplication by $e_1$. 
We will prove that
an analogous statement holds also for $J^{(k,r)}_n$. 

\begin{lem}\label{lem:3.1}
The ideal $J^{(k,r)}_n$ is invariant under the action of 
the operators $D_n(X;q,t)$, $E_m$ ($m\ge 0$). 
\end{lem}
\begin{proof}
We show that 
for any $f\in J^{(k,r)}_n$ and 
$I\subset\{1,\cdots,n\}$,   
$A_I(x;t)(T_{I}f)(x)$ satisfies the condition \eqref{zero}.
The assertion of the lemma is a corollary of this fact. 

Let $x_i$, $s_i$ be as in \eqref{zero}. 
Because of \eqref{qatb}, the denominator of $A_I(x;t)$ does not vanish. 
Set $\tilde{x}_i=qx_i$ if $i\in I$ and $\tilde{x}_i=x_i$ otherwise. 
Then $\tilde{x}_{\overline{i+1}}/\tilde{x}_i=tq^{\tilde{s}_i}$ for 
$1\le i\le k+1$, where 
$\overline{i}=i$ ($1\le i\le k+1$), $\overline{k+2}=1$, and
\be
\tilde{s}_i=\begin{cases}
s_i+1& (i\not \in I, \overline{i+1}\in I), \\
s_i-1& (i\in I, \overline{i+1}\not \in I), \\
s_i& (\mbox{ otherwise }). \\
\end{cases}
\en
If we have $s_i=0$ and 
$i\in I$, $\overline{i+1}\not\in I$
for some $1\le i\le k+1$, then 
$tx_i-x_{\overline{i+1}}=0$ and $A_I(x;t)=0$. 
Otherwise $\tilde{s}_i\ge 0$ for all $i$ and  
$\sum_{i=1}^{k+1}\tilde{s}_i=r-1$. 
Hence the wheel condition \eqref{zero} 
implies $(T_{I}f)(x)=f(\tilde{x})=0$.
\end{proof}

\subsection{Proof of an inclusion relation}\label{subsec:3.3}
In this subsection we prove the inclusion 
\bea
I^{(k,r)}_n\subset J^{(k,r)}_n. 
\label{incl}
\ena
If $n\le k$, then both sides 
are equal to $\Lambda_n$. 
Hence it suffices to consider the case $n\ge k+1$. 
\medskip

\noindent{\bf a) The case $n=k+1$}.\quad
{}For $i=1,\cdots,k+1$,  
let us call ($C_i$) the following statement: 
\be
\lambda_1-\lambda_{k+1}\ge r,
~~\ell(\lambda)\le i
\Longrightarrow 
P_\lambda\in J^{(k,r)}_{k+1}.
\en
We prove ($C_i$) by induction on $i$. 

In the case $i=1$, $\lambda$ has only one row.  
The corresponding Macdonald polynomials have 
a generating function given by a special case of 
the Cauchy identity (eq.(4.13), \cite{M}) 
\bea
\sum_{l\ge 0}P_{(l)}(x_1,\cdots,x_n;q,t)
\frac{(t;q)_l}{(q;q)_l}y^l
=\prod_{i=1}^n\frac{(tx_iy;q)_\infty}{(x_iy;q)_\infty},
\label{Cau}
\ena
where $(z;q)_m=\prod_{i=0}^{m-1}(1-zq^i)$. 
Let $s_1,\cdots,s_{k+1}\in\Z_{\ge 0}$ be integers satisfying 
$s_1+\cdots+s_{k+1}=r-1$. 
Let $x_i$, $s_i$ be as in \eqref{zero}, 
and specialize \eqref{Cau} accordingly, taking $n=k+1$. 
The right hand side becomes 
$
\prod_{i=1}^{k+1}(t^iq^{s_1+\cdots+s_{i-1}}x_1y;q)_{s_i} 
$,
which is a polynomial of degree $r-1$ in $y$. 
This implies that $P_{(l)}$ for $l\ge r$ 
vanishes under the condition \eqref{zero}.

Suppose that for some $i$ such that $2\le i\le k+1$, the condition
($C_{i-1}$) is true. We show ($C_{i}$) by induction on $\lambda_i>0$. 
Let $\lambda$ be as in ($C_{i}$) and 
set $\mu=\lambda_{(i)}$,  $\nu=\mu^{(i+1)}$. 
We have $P_\mu\in J^{(k,r)}_{k+1}$ by the induction hypothesis. 
{}From the formulas \eqref{e1}, \eqref{E2} and ($C_{i-1}$), 
we have modulo $J^{(k,r)}_{k+1}$ 
\be
e_1P_\mu
&\equiv &
\psi'_{\lambda/\mu}P_{\lambda}
+\psi'_{\nu/\mu}P_{\nu}, 
\\
c E_2P_\mu
&\equiv &
(1-q^{\mu_i}t^{1-i})\psi'_{\lambda/\mu}
P_{\lambda}
+
(1-t^{-i})\psi'_{\nu/\mu}
P_{\nu},
\en
where $c=t^{n-1}/(1-q)$. 
By Lemma \ref{lem:3.1}, the left hand sides belong to 
$J^{(k,r)}_{k+1}$. 
By Lemma \ref{lem:2.11} we have 
$\psi'_{\lambda/\mu}\neq 0$, and  
$1-q^{\mu_i}t^{1-i}\neq 1-t^{-i}$ by \eqref{qatb}. 
We conclude that $P_{\lambda}$ belongs to $J^{(k,r)}_{k+1}$. 
\medskip

\noindent{\bf b) The case $n\ge k+2$}.\quad
Let 
\be
\rho:\Lambda_n\rightarrow \Lambda_{n-1},\quad 
\rho(f)(x_1,\cdots,x_{n-1})=f(x_1,\cdots,x_{n-1},0) 
\en
denote the specialization map. 
To prove \eqref{incl} for $n\ge k+2$, it suffices to show that 
\bea
f\in I^{(k,r)}_n\Longrightarrow 
\rho(\partial_n^j f)\in J^{(k,r)}_{n-1}
\qquad (j\ge 0), 
\label{rho}
\ena
where $\partial_n=\partial/\partial x_n$.

Since $P_\lambda(x_1,\cdots,x_{n-1},0;q,t)
=P_\lambda(x_1,\cdots,x_{n-1};q,t)$, 
\eqref{rho} holds for $j=0$. 
Suppose it is true for $j-1$. 
Set 
\bea
\partial_n^{j-1}E_0f=\sum_{i=1}^{n}X_i,
\label{Xi}
\ena
where 
\be
&&X_i=A'_i\partial_n^{j-1}
\left(\frac{tx_i-x_n}{x_i-x_n}
\frac{\partial f}{\partial_q x_i}\right),
\qquad A'_i=\prod_{l(\neq i,n)}\frac{tx_i-x_l}{x_i-x_l},
\\
&&X_n=
\partial_n^{j-1}
\left(A_n\frac{\partial f}{\partial_q x_n}\right).
\en
Since $E_0I^{(k,r)}_n\subset I^{(k,r)}_n$, 
by induction hypothesis the image by $\rho$ of 
the left hand side of \eqref{Xi} 
belongs to $J^{(k,r)}_{n-1}$. 

Consider the terms with $i\le n-1$, 
\be
\rho(X_i)&=&
\sum_{s=0}^{j-1}\binom{j-1}{s}
\rho\Bigl(\partial^{j-1-s}_n\frac{tx_i-x_n}{x_i-x_n}\Bigr)
A_i'\frac{\partial}{\partial_q x_i}\rho(\partial_n^s f).
\en
Since $g=\rho(\partial_n^s f)$ belongs to $J_{n-1}^{(k,r)}$, 
$A_i'\frac{\partial}{\partial_q x_i}g$ satisfies \eqref{zero}, 
as we have seen in the proof of Lemma \ref{lem:3.1}. 
This shows that $\rho(X_i)\in J^{(k,r)}_{n-1}$ ($i\le n-1$). 
Noting that
\be
\rho\Bigl(\partial_n^{s-1}\frac{\partial f}{\partial_q x_n}\Bigr)
=a_s\rho(\partial_n^{s}f)
\en
with $a_{s}=(1-q^s)/(s(1-q))$, we find
\be
\rho(X_n)&=&\sum_{s=0}^{j-1}\binom{j-1}{s}
\rho(\partial_n^{j-1-s}A_n)
\rho\Bigl(\partial_n^{s}\frac{\partial f}{\partial_q x_n}\Bigr)
\\
&\equiv& a_j\rho(\partial_n^j f)
~~\bmod~~ J^{(k,r)}_{n-1}.
\en
The assertion \eqref{rho} follows from these. 

\subsection{Comparison of dimensions}\label{subsec:3.4}
In this subsection we finish the proof of 
Theorem \eqref{thm:1} by comparing dimensions. 

Counting $\mathop{\rm deg}x_i=1$ for all $i$, 
we denote by $\Lambda_{n,d}\subset \Lambda_n$ the 
subspace of polynomials of homogeneous degree $d$.  
The spaces 
$J^{(k,r)} $, $I^{(k,r)}$ are homogeneous with respect to 
the bi-grading of $\Lambda=\oplus\Lambda_{n,d}$. 
We set $J^{(k,r)}_{n,d}=J^{(k,r)}\cap \Lambda_{n,d}$, 
$I^{(k,r)}_{n,d}=I^{(k,r)}\cap \Lambda_{n,d}$.  
Let also $\pi_{n,d}$ be the set of partitions 
$\lambda\in\pi_n$ satisfying 
$\sum_{i=1}^{n}\lambda_i=d$. 
We denote by $\mathcal{C}^{(k,r)}_{n}$ the set of 
$(k,r,n)$-admissible partitions,  
and $\mathcal{C}^{(k,r)}_{n,d}=\mathcal{C}^{(k,r)}_{n}
\cap \pi_{n,d}$.

Let us specialize the parameters further to $t=1,q=\tau$,  
where $\tau$ is a primitive $(r-1)$-th root of unity.
Under this specialization, the wheel condition \eqref{zero} becomes 
\bea
\mbox{$f=0$ if $x_{i}=\tau^{p_i}x_1$ 
for all $p_i\in\Z$\quad ($2\le i\le k+1$).} 
\label{dzero}
\ena
We consider the corresponding polynomial space 
$\overline{J}^{(k,r)}=\oplus_{n\ge 0}\overline{J}^{(k,r)}_n$, where
\be
\overline{J}^{(k,r)}_n
=\{ f\in\C[x_1,\cdots,x_n]^{S_n}
\mid \mbox{ $f$ satisfies \eqref{dzero} }\}. 
\en

Let us determine the character of $\overline{J}^{(k,r)}$.
{}For this purpose it is convenient to pass to the dual space.
Let $\overline{R}=\C[e_0,e_1,\cdots]$ be the polynomial ring in 
indeterminates $\{e_i\}_{i\ge 0}$, 
equipped with the bi-grading $\mathop{{\rm deg}}e_i=(1,i)$. 
Let 
$e(\zeta)=\sum_{i\ge 0}e_i\zeta^i$ be the generating series.  
Consider the ideal $\overline{\mathcal{J}}\subset \overline{R}$ generated by the Fourier 
coefficients of 
\bea
e(\zeta)e(\tau^{p_2}\zeta)\cdots e(\tau^{p_{k+1}}\zeta),
\label{aaa}
\ena
where $p_2,\cdots,p_{k+1}$ run through arbitrary integers. 
A standard argument  
shows (see e.g. \cite{FJMMT}, Lemma 3.3) 
that there is a non-degenerate coupling 
\be
(\overline{R}/\overline{\mathcal{J}})\times \overline{J}^{(k,r)} \rightarrow \C, 
\en
through which each
homogeneous component $(\overline{R}/\overline{\mathcal{J}})_{n,d}$ 
is isomorphic to $\bigl(\overline{J}^{(k,r)}_{n,d}\bigr)^*$.

\newcommand{\bnu}{{\boldsymbol \nu}}
In the below, we write 
\be
e_\lambda=e_{\lambda_1}e_{\lambda_2}\cdots e_{\lambda_n} 
\en
for a partition $\lambda=(\lambda_1,\cdots,\lambda_n)\in\pi_n$.
 
\begin{prop}\label{prop:3.3}
The space $\overline{R}/\overline{\mathcal{J}}$ is spanned by the set of monomials 
\bea
B=\{e_{\lambda}\mid \lambda\in \mathcal{C}^{(k,r)}_n\}.
\label{mon}
\ena
\end{prop}
\begin{proof}
We use the lexicographic ordering 
$\mu\succ\lambda$ defined by 
$\mu_1=\lambda_1,\cdots,\mu_{i-1}=\lambda_{i-1},\mu_i>\lambda_i$ for some $i$. 

If we write 
\be
e(\zeta)=\sum_{j=0}^{r-2}\zeta^j e_{1\,j+2}(\zeta^{r-1}),
\quad
e_{1\,j+2}(z)=\sum_{l\ge 0}e_{(r-1)l+j} z^l, 
\en
then the defining relations \eqref{aaa} 
for the ideal $\overline{\mathcal{J}}$ can be stated equivalently as 
\bea
e_{12}(z)^{\nu_0}e_{13}(z)^{\nu_1}\cdots
e_{1r}(z)^{\nu_{r-2}}\equiv 0
\quad \bmod \overline{\mathcal{J}}.
\label{int2}
\ena
Here $\bnu=(\nu_0,\cdots,\nu_{r-2})$ runs over 
non-negative integers $\nu_j\in\Z_{\ge 0}$ such that 
$\sum_{j=0}^{r-2}\nu_j=k+1$.
Hence $\overline{\mathcal{J}}$ is generated by elements of the form 
$\sum_{\mu\in K_d(\bnu)}C_\mu e_\mu$,
where $C_\mu$ are positive integers and 
\be
K_d(\bnu)=\{\mu\in\pi_{k+1,d}\mid 
\sharp\{i\mid\mu_i\equiv a \bmod r-1\}=\nu_a
~~(0\le a\le r-2)
\}.
\en
It is easy to see that, for each $d$ and $\bnu$, 
the set $K_d(\bnu)$ is either empty or 
contains a unique element
$\lambda$ satisfying $\lambda_1-\lambda_{k+1}\le r-1$.  
If $\mu\in K_d(\bnu)$ with $\mu\neq \lambda$, then 
we have $\mu\succ \lambda$.
Moreover, any ($k,r,k+1$)-non-admissible partition 
$\lambda$ is the minimal element of some $K_d(\bnu)$.
Therefore, for such $\lambda$, 
$e_\lambda$ belongs to the linear span of $B$. 

Let us consider the general case $\lambda\in\pi_{n,d}$. 
Suppose $\lambda_i-\lambda_{i+k}\le r-1$ for some $i$. 
Set $\mu=(\lambda_i,\cdots,\lambda_{i+k})$, and 
let $\tilde{\lambda}$ be the partition obtained from $\lambda$ by deleting $\mu$. 
We have 
\be
e_\lambda=e_{\tilde{\lambda}}e_\mu.
\en
Rewriting $e_\mu$ as a linear span of elements of $B$, we obtain 
\be
e_\lambda=\sum_{\nu\in\pi_{n,d}\atop \nu\succ\lambda}
c_{\lambda\nu}e_\nu
\qquad (c_{\lambda\mu}\in\C). 
\en
If non-admissible $\nu$ appears in the right hand side, then 
we can repeat the same procedure for $\nu$.
Since $\pi_{n,d}$ is a finite set, 
this process terminates after a finite number of steps, 
giving $e_\lambda$ as a linear span of $B$. 
\end{proof}

Let us return to the proof of Theorem \ref{thm:1}. 

{}From the definition of $I^{(k,r)}_n$, it is clear that
\bea
\sharp \mathcal{C}^{(k,r)}_{n,d}
=\dim_K I^{(k,r)}_{n,d}.
\label{Idim}
\ena
The inclusion \eqref{incl} implies 
\bea
\dim_K I^{(k,r)}_{n,d}
\le \dim_K J^{(k,r)}_{n,d}. 
\label{dim}
\ena
Since the defining relations \eqref{dzero} 
are obtained by specializing \eqref{zero},  
we have 
\bea
\dim_K J^{(k,r)}_{n,d}
\le \dim_\C \overline{J}^{(k,r)}_{n,d}.
\label{est1}
\ena
{}From Proposition \ref{prop:3.3} we find 
\bea
\dim_\C \overline{J}^{(k,r)}_{n,d}
=\dim_\C(\overline{R}/\overline{\mathcal{J}})_{n,d}
\le\sharp \mathcal{C}^{(k,r)}_{n,d}.
\label{est2}
\ena
Combining \eqref{Idim}--\eqref{est2}
we conclude that the equality takes place in \eqref{dim}. 
Proof of Theorem \ref{thm:1} is now complete. 

\begin{cor}\label{cor:3.1}
The set $B$ is a basis of $\overline{R}/\overline{\mathcal{J}}$. 
We have 
\be
\dim_K J^{(k,r)}_{n,d}
=\dim_\C \overline{J}^{(k,r)}_{n,d}
=\sharp \mathcal{C}^{(k,r)}_{n,d}.
\en
\end{cor}
\medskip

\section{Monomial basis}
In \cite{P}, a monomial base of the form 
\eqref{mon} was constructed for arbitrary irreducible 
integrable representations of $\slrh$.  
As we mentioned in Introduction, the space $\overline{R}/\overline{\mathcal{J}}$ 
is a principal subspace of an integrable representation. 
We have given, for this special case, 
an alternative proof of the result of \cite{P} 
and its deformation at the same time.   
In this Section we present 
a counterpart of the monomial basis for 
general integrable representations of $\slrh$.  

Let $R=K[\{e_i\}_{i\ge 0}]$. 
{}For an array $(b_0,\cdots,b_{r-2})$ 
of non-negative integers satisfying 
$0\le b_0\le\cdots\le b_{r-2}\le k$, 
introduce the quotient space 
\be
W_{b_0,b_1,\cdots,b_{r-2}}=R/I_{b_0,b_1,\cdots,b_{r-2}}, 
\en
where $I_{b_0,b_1,\cdots,b_{r-2}}$ denotes the ideal of $R$ 
generated by 
\bea
&&
\mbox{the Fourier coefficients of 
$e(z)e(tq^{s_1}z)\cdots e(tq^{s_k}z)$}
\quad (s_i\ge 0,\sum_{i=1}^ks_i\le r-1),
\label{ideal1}\\
&&
e_0^{a_0}e_1^{a_1}\cdots e_{r-2}^{a_{r-2}},
\quad
\mbox{$a_0+\cdots+a_i>b_i$ for some $0\le i\le r-2$}.
\label{ideal2}
\ena
The space $W_{b_0,b_1,\cdots,b_{r-2}}$ is bi-graded. 
We set
\be
&&\chib_{b_0,b_1,\cdots,b_{r-2}}(v,z)
=\sum_{n,d}\dim \left(W_{b_0,b_1,\cdots,b_{r-2}}(v,z)\right)_{n,d}v^dz^n.
\en

In order to parameterize the monomial basis in general, 
we find it convenient to represent a partition $\lambda$ 
by the numbers $a_i$ of parts $i$ of $\lambda$, 
where we set $a_i=0$ for $i>\lambda_1$.  
We have a sequence of non-negative integers $a=(a_i)_{i=0}^{\infty}$ with 
$a_i=0$ for $i$ large enough.
Let $\mathcal{C}_{b_0,b_1,\cdots,b_{r-2}}$ 
denote the set of all such sequences satisfying the following conditions.
\be
&&a_i+a_{i+1}+\cdots+a_{i+r-1}\le k~~(i\ge 0),
\\
&&
a_0\le b_0,~~
a_0+a_1\le b_1,\cdots a_0+\cdots+a_{r-2}\le b_{r-2}.
\en
Define its character by
\be
\chiC_{b_0,b_1,\cdots,b_{r-2}}(v,z)=
\sum_{a\in \mathcal{C}_{b_0,b_1,\cdots,b_{r-2}}}
v^{\sum_{i\ge 0}ia_i}z^{\sum_{i\ge 0}a_i}.
\en

Our aim is to show the following. 
\begin{prop}\label{prop:302}
The space $W_{b_0,b_1,\cdots,b_{r-2}}$ has a monomial basis 
\bea
B_{b_0,b_1,\cdots,b_{r-2}}=\{
\prod_{i=0}^\infty e_i^{a_i}
\mid 
a=(a_i)_{i=0}^\infty
\in 
\mathcal{C}_{b_0,b_1,\cdots,b_{r-2}}\}.
\label{mon302}
\ena
\end{prop}

{}From Corollary \eqref{cor:3.1} we know that 
$W_{k,k,\cdots,k}$ has the set \eqref{mon} as monomial basis.   
Therefore the space 
$W_{b_0,b_1,\cdots,b_{r-2}}$ is spanned by the monomials
\eqref{mon302}, and we have
\bea
&&
\chib_{b_0,b_1,\cdots,b_{r-2}}(v,z)
\le
\chiC_{b_0,b_1,\cdots,b_{r-2}}(v,z).
\label{eq1}
\ena
Here and in what follows, 
for two formal series $f=\sum_{n,d}f_{n,d}q^dz^n$ and 
$g=\sum_{n,d}g_{n,d}q^dz^n$, we write 
$f\le g$ to mean $f_{n,d}\le g_{n,d}$ for all $n,d$. 

The following Lemma is immediate.
\begin{lem}\label{lem:11}
We have the recursion relation
\bea
&&\chiC_{b_0,b_1,\cdots,b_{r-2}}(v,z)
\nn\\
&&=
\chiC_{b_0-1,b_1,\cdots,b_{r-2}}(v,z)
+
z^{b_0}\chiC_{b_1-b_0,\cdots,b_{r-2}-b_0,k-b_0}(v,vz).
\label{recur}
\ena
\end{lem}

\begin{lem}\label{lem:1}
We have an exact sequence
\be
W_{b_1-b_0,\cdots,b_{r-2}-b_0,k-b_0}
\overset{\varphi}{\longrightarrow}
W_{b_0,b_1,\cdots,b_{r-2}}
\overset{\pi}{\longrightarrow}
W_{b_0-1,b_1,\cdots,b_{r-2}}
\longrightarrow
0.
\en
Here $\pi$ is the canonical surjection and 
$\varphi(u)=e_0^{b_0} T(u)$,
where $T:R\rightarrow R$ 
is a homomorphism of algebras defined by $T(e_i)=e_{i+1}$. 
\end{lem}
\begin{proof}
Let us show that $\varphi$ is well defined. 
Abusing the notation we use the same letter for the map 
$\varphi:R\rightarrow R$ defined by the formula above.
Noting that $T(e(z))=z^{-1}(e(z)-e_0)$, we have
\be
\varphi\bigl(\prod_{i=0}^ke(tq^{s_i}z)\bigr)
&=&z^{-k-1}e_0^{b_0}\prod_{i=0}^k\bigl(e(tq^{s_i}z)-e_0\bigr)
\\
&\equiv&z^{-k-1}e_0^{b_0}\prod_{i=0}^ke(tq^{s_i}z)
\\
&\equiv&0\quad
\bmod I_{b_0,\cdots,b_{r-2}}.
\en
Similarly, if $a_0+\cdots+a_{i-1}>b_i-b_0$ for some $i$ 
(where we set $b_{r-1}=k$), then
\be
\varphi\bigl(\prod_{i=0}^{r-2}e_i^{a_i}\bigr)
&=&e_0^{b_0}\bigl(\prod_{i=1}^{r-1}e_i^{a_{i-1}}\bigr)
\\
&\equiv&0\quad
\bmod I_{b_0,\cdots,b_{r-2}}. 
\en
This shows that 
$\varphi(I_{b_1-b_0,\cdots,b_{r-2}-b_0,k-b_0})\subset 
I_{b_0,\cdots,b_{r-2}}$. 

Clearly $I_{b_0-1,\cdots,b_{r-2}}/I_{b_0,\cdots,b_{r-2}}$ is spanned by 
$e_0^{b_0}\prod_{i\ge 1}e_i^{a_i}$ with 
$b_0+a_1+\cdots+a_i\le b_i$ for $1\le i\le r-1$.
The exactness follows from this. 
\end{proof}

Proposition \ref{prop:302} is a consequence of the following. 
\begin{prop}
We have
\be
\chib_{b_0,b_1,\cdots,b_{r-2}}(v,z)
=
\chiC_{b_0,b_1,\cdots,b_{r-2}}(v,z).
\en
\end{prop}
\begin{proof}
Lemma \ref{lem:1} implies that 
\be
&&\chib_{b_0,b_1,\cdots,b_{r-2}}(v,z)
\nn\\
&&\le
\chib_{b_0-1,b_1,\cdots,b_{r-2}}(v,z)
+
z^{b_0}\chib_{b_1-b_0,\cdots,b_{r-2}-b_0,k-b_0}(v,vz).
\en
Taking $b_0=l,b_1=\cdots=b_{r-2}=k$ we find for $0<l\le k$ that 
\be
\chib_{l,k,\cdots,k}(v,z)
&\le &
\chib_{l-1,k,\cdots,k}(v,z)
+
z^l\chib_{k-l,k-l,\cdots,k-l}(v,vz)
\\
&\le &
\chiC_{l-1,k,\cdots,k}(v,z)
+
z^l\chiC_{k-l,k-l,\cdots,k-l}(v,vz)
\\
&=&\chiC_{l,k,\cdots,k}(v,z).
\en
In the last line we used \eqref{recur}.
{}From Corollary \ref{cor:3.1} we have 
\bea
\chib_{k,k,\cdots,k}(v,z)=
\chiC_{k,k,\cdots,k}(v,z).
\label{eq3}
\ena
Using \eqref{eq3} as a base of induction on $l=k,k-1,\cdots,0$,  
we obtain
\be
&&\chib_{l,k,\cdots,k}(v,z)=\chiC_{l,k,\cdots,k}(v,z),
\\
&&
\chib_{k-l,\cdots,k-l}(v,z)=\chiC_{k-l,\cdots,k-l}(v,z).
\en

Arguing similarly, we find by induction that 
\be
\chib_{
b,\cdots,b,
b_{s},\cdots,b_{r-2}}(v,z)
=
\chiC_{
b,\cdots,b,
b_{s},\cdots,b_{r-2}}(v,z)
\en
for all
$1\le s\le r-2$ and $0\le b\le b_s\le\cdots\le b_{r-2}\le k$. 
\end{proof}

\newcommand{\W}{\mathcal{W}}

\section{Discussions}\label{sec:4}
In this Section we discuss a possible connection 
between this paper and the representations of 
the $\W_k$ algebra associated with $\slkh$. 

Recall the following well-known phenomenon 
in representation theory of $\slth$. 
Let $L$ be an irreducible integrable representation 
of $\slth$ of level $k$. 
Let $\hh$ be the Heisenberg subalgebra of $\slth$. 
Then we have a decomposition 
$L=\oplus_{\alpha\in\Z} \pi_\alpha\otimes S_\alpha$, 
where $\pi_\alpha$ are irreducible representations of 
$\hh$ and $S_\alpha$ are irreducible representations 
of the $\W_k$ algebra in the minimal series ($k+1,k+2$). 
In other words, we have on $L$ 
an action of $\W_k$ commuting with $\hh$. 
Therefore, in some sense $\slth$ on level $k$ is an
``extension'' of the product $\W_k\times\hh$. 

We suggest a possible generalization of this construction. 
Consider a $\W_k$ minimal series representation 
labeled by relatively prime integers $(p,q)$. 
We will write them in the form $(p,q)=(k+s,k+r)$. 
There is a special set of primary fields 
$\varphi_0(z)={\rm id},\varphi_1(z),\cdots,\varphi_{k-1}(z)$ 
with the operator product expansion 
\be
\varphi_\alpha(z)\varphi_\beta(w)=
(z-w)^{\Delta_{\alpha\beta}}\varphi_\gamma(w)+\cdots.
\en
Here $\gamma\equiv \alpha+\beta\bmod k$, 
$\Delta_{\alpha\beta}$ 
are some rational numbers,  
and the higher terms denoted by the dots involve the descendants of 
$\varphi_\gamma(w)$. 
The operators $\{\varphi_\alpha(z)\}$
constitute a generalization of the parafermion algebra \cite{LW}. 
Now consider a one-component Heisenberg algebra $\hh$ 
which commutes with $\W_k$. 
Let $V_+(z),V_-(z)$ be vertex operators for $\hh$ 
with the properties 
\be
&&
V_\pm(z)V_\pm(w)=(z-w)^{-\Delta_{11}}:V_\pm(z)V_\pm(w):,
\\
&&
V_\pm(z)V_\mp(w)=(z-w)^{\Delta_{11}}:V_\pm(z)V_\mp(w):.
\en
Set 
\be
&&e(z)=V_+(z)\varphi_1(z), 
\\
&&f(z)=V_-(z)\varphi_{k-1}(z).
\en
It is easy to see that 
$[e(z),e(w)]=0$, $[f(z),f(w)]=0$, 
and that $[e(z),f(w)]=0$ for $z\neq w$. 

Let $\A^k_{k+s,k+r}$ be the vertex operator algebra generated by
$e(z),f(z)$. 
We view $s,r$ as parameters of the algebra 
and $k$ as the level of the representation.
{}For instance, the algebra $\A^k_{k+1,k+2}$ is $\slth$ acting on 
integrable representations of level $k$.  
Not much is known about this algebra in general.
%

When $k$ is a non-negative integer, 
we can impose some additional integrability conditions. 
The integrability conditions can be reformulated as 
a statement about the matrix elements of the current $e(z)$. 

{}For example, consider a 
level $k$ integrable representation $L$ of $\slth$.
The matrix elements of $e(z)$, 
\bea
P(z_1,\cdots,z_n)=\langle v^{\vee},e(z_1)\cdots e(z_n) v\rangle
\qquad (v\in L, v^\vee\in L^*), 
\label{matel}
\ena
are symmetric Laurent polynomials in $(z_1,\cdots,z_n)$ 
satisfying
\bea
\mbox{$P=0$ if $z_1=\cdots=z_{k+1}$}. 
\label{kdiag}
\ena
The integrable irreducible representations 
can be realized in a space of symmetric Laurent polynomials 
in infinite set of variables satisfying the zero 
condition \eqref{kdiag} \cite{FS}. 
Therefore, for $\slth$, 
we can start from the relation $e(z)^{k+1}=0$ 
and reconstruct the level $k$ vacuum representation. 
It is possible further to 
find the structure of a vertex operator algebra on it.

{}For general $\A^k_{k+s,k+r}$,  
the matrix elements are also symmetric Laurent polynomials 
with some zero conditions on the diagonal of codimension $k+1$. 
As in the case of $\slth$, one expects that these 
zero conditions determine the integrable representations, 
and even the algebra $\A^k_{k+s,k+r}$ itself. 
When we try to establish these facts, 
the first obstacle is that 
the integrability conditions are not known. 
In the special case $\A^k_{k+1,k+r}$, 
we have presented in \cite{FJMM} 
a conjectural description of the integrability 
condition in the language of matrix elements. 
More precisely, 
we considered certain subspaces (in fact, ideals) 
of symmetric Laurent polynomials 
which are spanned by Jack polynomials. 
We expect that they coincide with the space 
of all matrix elements of the current $e(z)$ 
in integrable representations of $\A^k_{k+1,k+r}$.

In this paper, we found some evidence that all 
this can and should be `$q$-deformed'. 
An obvious counterpart of $\W_k$ would be
its elliptic deformation. 
We conjecture that the algebra $\A^k_{k+s,k+r}$ 
can be $q$-deformed, 
in such a way that the currents $e(z),f(z)$ remain commutative:
$[e(z),e(w)]=0$, $[f(z),f(w)]=0$. 
In the undeformed case, the matrix elements of 
$e(z)$ and $f(z)$ have the structure
\be
\langle w^\vee, e(z_1)f(z_2)w \rangle
=F(z_1,z_2)(z_1-z_2)^{-m},
\en
where $F$ is a Laurent polynomial.
If $s=1$, then $m=r$. 
After deformation, we expect that they take the form
$F(z_1,z_2)\prod_{i=1}^m(z_1-q_iz_2)^{-1}$
with some $q_1,\cdots,q_m$. 

In this paper
we have found an explicit description of 
the integrability condition for the current $e(z)$ itself in 
the special case $s=1$.
The analogue description is only implicit in the undeformed case, 
see \cite{FJMM}.
This is one of the advantages of considering the $q$-deformation. 

Note that for a special value of the parameters $(q,t)=(\tau,1)$, 
the algebra $\A^k_{k+1,k+r}$ is known -
somewhat surprisingly, 
it is $\widehat{\mathfrak{sl}}_r$ at level $k$. 
\bigskip

\noindent
{\it Acknowledgments.}\quad 
B.F. is partially supported by the grants,
CRDF RP1-2254, INTAS 00-55, RFBR 00-15-96579.
M.J. is partially supported by
the Grant-in-Aid for Scientific Research (B2)
no.14340040, and 
T.M. is partially supported by
the Grant-in-Aid for Scientific Research (A)
no. 13304010, Japan Society for the Promotion of Science.
E.M. is partially supported by NSF grant DMS-0140460.



\begin{thebibliography}{FJMMT}

\bibitem[FJLMM]{FJLMM}
B.~Feigin, M.~Jimbo, S.~Loktev, T.~Miwa and E.~Mukhin, 
{Bosonic formulas for $(k,l)$--admissible partitions},
math.QA/0107054, to appear in Ramanujan J, 2001.

\bibitem[FJMM]{FJMM}
B.~Feigin, M.~Jimbo, T.~Miwa and E.~Mukhin, 
{\it A differential ideal of symmetric polynomials 
spanned by Jack polynomials at $\beta=-(r-1)/(k+1)$}, 
Int. Math. Res. Notice. {\bf 23} 1223--1237 (2002).

\bibitem[FJMMT]{FJMMT}
B.~Feigin, M.~Jimbo, T.~Miwa, E.~Mukhin and Y.~Takeyama,
{\it Symmetric polynomials vanishing on the diagonals shifted by roots of unity},
in preparation. 

\bibitem[FO]{OF}
B.~L. Feigin and A.~Odesskii,
\newblock
Vector bundles on elliptic curve and Sklyanin algebras
\newblock
Amer. Math. Soc. Transl. Ser. 2, {\bf 185} (1998), 65--84.

\bibitem[FS]{FS}
B.~L. Feigin and A.~V. Stoyanovsky,
\newblock Functional models 
for representations of current algebras and 
  semi-infinite {Schubert} cells, 
\newblock {\em Funct.~Anal.~and~Its~Appl.}
{\bf 28} (1993) 55--72.

\bibitem[Int]{Int} 
A.~Kirillov, M.~Noumi, {\it Affine Hecke algebras
and raising operators for Macdonald polynomials}, 
Duke Math.J. {\bf 93} (1998) 1--39. \\
S.~Sahi, {\it Interpolation, integrality, and a generalization of
Macdonald's polynomials},
Int. Math. Res. Notices {\bf 10} (1996) 457-471.\\
K.~Friedrich, {\it Integrality of two variable Kostka functions}, 
J. Reine Angew. Math. {\bf 482} (1997) 177--189.\\
A.~Garsia and J.~Remmel, {\it Plethystic formulas and positivity for
$q,t$-Kostka coefficients},
Progr. Math., Birkh\"{a}user, {\bf 161} (1998) 245--262.\\
A.~Garsia and G.~Tesler, {\it Plethystic formulas for
$q,t$-Kostka coefficients},
Adv. Math. {\bf 123} (1996) 144--222.  
  
\bibitem[L]{L} M.~Lassalle, 
{\it Coefficients binomiaux g\'{e}n\'{e}ralis\'{e}s et polyn\^{o}mes de Macdonald},
J. Funct. Anal. {\bf 158} (1998) 289--324.

\bibitem[LW]{LW}
J.~Lepowsky and R.~L. Wilson, 
\newblock Construction of the affine lie algebra {$A^{(1)}_1$}, 
\newblock {\em Commun.~Math.~Phys.} {\bf 62} (1978) 43--53. 

\bibitem[M]{M}I.~Macdonald, 
{\it Symmetric functions and Hall polynomials, 2nd ed.}, 
Oxford University Press, New York, 1995.

\bibitem[P]{P} M.~Primc, 
{\it Vertex operator construction of standard modules for $A^{(1)}_n$}, 
Pacific J. Math. {\bf 162} (1994) 143--187.

\end{thebibliography}
\end{document}